\documentclass{siamart1116}

\usepackage{amsfonts}
\usepackage{graphicx}
\usepackage{pgfplots}
\usepackage{tabu}
\ifpdf
  \DeclareGraphicsExtensions{.eps,.pdf,.png,.jpg}
\else
  \DeclareGraphicsExtensions{.eps}
\fi

\numberwithin{theorem}{section}

\newcommand{\TheTitle}{A High-Order Lower-Triangular Pseudo-Mass Matrix for Explicit Time Advancement of \texorpdfstring{$hp$}{hp} Triangular Finite Element Methods.} 
\newcommand{\RunningTitle}{A High-Order Lower-Triangular Pseudo-Mass Matrix}
\newcommand{\TheAuthors}{Jay M. Appleton and B. T. Helenbrook}
\newcommand{\TheSubDate}{06/25/2019}

\headers{\RunningTitle}{\TheAuthors}

\title{{\TheTitle}\thanks{Submitted to SIAM Journal on Numerical Analysis \TheSubDate.}}

\author{
  Jay Miles Appleton\thanks{Department of Mathematics, Clarkson University, 8 Clarkson Ave. Box 5815 Potsdam, NY
    (\email{appletj@clarkson.edu}).}
  \and
  Brian T. Helenbrook\thanks{Department of Mechanical and Aeronautical Engineering, Clarkson University, 8 Clarkson Ave Box 5725 
    Potsdam, NY (\email{helenbrk@clarkson.edu}, \url{https://www.clarkson.edu/mae/faculty_pages/helenbrook.html}).}
}

\usepackage{amsopn}

\ifpdf
\hypersetup{
  pdftitle={\TheTitle},
  pdfauthor={\TheAuthors}
}
\fi

\begin{document}
\maketitle

  \begin{abstract} 
    Explicit time advancement for continuous finite elements requires the inversion of a global mass matrix.  For spectral element simulations on quadrilaterals and hexahedra, there is an accurate approximate mass matrix which is diagonal, making it computationally efficient for explicit simulations.   In this article it is shown that for the standard space of polynomials used with triangular elements, denoted $\mathcal{T}(p)$ where $p$ is the degree of the space, there is no diagonal approximate mass matrix that permits accurate solutions.  Accuracy is defined as giving an exact projection of functions in $\mathcal{T}(p-1)$.   In light of this, a lower-triangular pseudo-mass matrix method is introduced and demonstrated for the space $\mathcal{T}(3)$.
    The pseudo-mass matrix and accompanying high-order basis allow for computationally efficient time-stepping techniques without sacrificing the accuracy of the spatial approximation for unstructured triangular meshes.
  \end{abstract}


  \begin{AMS}
    {65M60, 65D30, 65M20, 65N35, 65N40}
  \end{AMS}

 \begin{section}{Introduction}\label{Sec:Introduction}
 
    The two-dimensional spectral element method (SEM) discussed in \cite{Karniadakis2005} is a continuous Galerkin finite element method (FEM) for quadrilateral meshes that utilizes a high degree polynomial basis.  It is often coupled with explicit finite differences in time for the simulation of unsteady problems.  Explicit simulations are computationally efficient for many problems and also easy to implement in parallel \cite{Sun2000, Mavriplis1997}.  For FEM, explicit time advancement schemes typically require the inversion of a global mass matrix at each time step, which can be computationally expensive \cite{Abgrall2017b}.  This expense is avoided in the SEM because it uses Gauss-Lobatto (GL) integration, which provides the numerical integration scheme, the nodes (the GL points), and the basis functions (Lagrange interpolants over the GL points) for the method \cite{Pasquetti2017}.  This approach constitutes a nodal collocation FEM because each basis function is only nonzero at one GL integration point.  When the SEM mass matrix is evaluated using  GL numerical integration, it becomes diagonal, which eliminates the expense of the mass matrix inversion.  Note that with exact integration, the mass matrix is not diagonal.  The diagonality is a consequence of the numerical integration; there is no polynomial basis with continuity properties appropriate for continuous FEM that is exactly orthogonal (i.e. a basis with a diagonal mass matrix).  The SEM approximate mass matrix, which we call a pseudo-mass matrix, has the property that when used to project a function into the SEM approximation space it is exact if the function to be projected is a polynomial in the SEM space of one less degree than the approximation space \cite{Helenbrook2009}.   
 
    Although the SEM has become a popular numerical technique for simulating challenging problems in complicated domains \cite{Deville2002, Canuto2007}, the need for high-order methods on unstructured meshes with robust adaptivity has motivated the development of a triangular finite element method (TFEM) that can compare to the SEM \cite{Samson2012, Abgrall2017a, Hesthaven2000}.  A comparable TFEM would allow an implementation of numerical methods for differential equations that possesses high-order spatial accuracy, low-memory usage, and efficient parallelizability.  Such a method would lend itself well to structural mechanics \cite{Karniadakis2005}, direct numerical simulations of computational fluid dynamics \cite{Karniadakis2005, Helenbrook2009, Sherwin1995, Chen2011}, atmospheric modeling \cite{Giraldo2005, Eskilsson2004}, etc. \cite{Sousa2016, Zarei2016, Nicolas1999, Pasquetti2017}.  This has driven research into high-order continuous Galerkin TFEM possessing diagonal mass matrices that maintain a high level of accuracy \cite{Giraldo2006, Brazell2013, Taylor2000a, Cohen2001, Hesthaven1998}.  In \cite{Helenbrook2009}, Helenbrook showed that there is no nodal $C^0$ TFEM comparable to the SEM; specifically, he showed that there is no GL integration rule for triangles that can be used to create a $p^{th}$ degree nodal basis and is also accurate to order $2p-1$, which is the case on quadrilaterals.  An artifact of the proof was the derivation a unique set of modal vertex functions that provide an accurate diagonal vertex block of a pseudo-mass matrix.  The existence of these modal vertex functions inspire this investigation into an entire modal basis with a diagonal pseudo-mass matrix.   
  
    To investigate whether this basis and diagonal pseudo-mass matrix pair exists, the paper begins with a derivation of an explicit $C^0$ TFEM for an arbitrary basis and introduces the continuity requirements for a high-order basis.  The modified Dubiner basis \cite{Dubiner1991} is presented as a popular choice that satisfies these requirements.  The concept of a pseudo-mass matrix and its accuracy requirements are then mathematically defined.  A change of basis is introduced that can be used to map the Dubiner basis to any other high-order basis suitable for $C^0$ TFEM thus allowing the accuracy requirements to be defined for all bases.  Then, these requirements are used to show that on triangles there is no basis and diagonal pseudo-mass matrix pair that satisfies the accuracy constraints.  

To provide an alternative approach, we relax the diagonality requirement and instead look for a pseudo-mass matrix that can be inverted with only local operations.  To this end, a new lower-triangular pseudo-mass matrix for $p=3$ is introduced that provides a desirable alternative to the full mass matrix approach by avoiding the need for a global mass matrix.  This new method not only serves as a demonstration of concept, but is both a viable higher-order option and a segue into future work toward arbitrarily high-order continuous methods for triangles.
 \end{section}

\begin{section}{A High-Order Continuous Explicit TFEM}\label{Sec:Background_Stuff}
    This section provides the formulation of a high-order $C^0$ TFEM for a transient problem, thereby introducing the necessary concepts and definitions for the following sections.  The primary concern of an explicit time advancement method is the handling of the transient term.  Of particular interest in this paper are transient partial differential equations (PDEs) that can be separated into a spatial operator ${\cal L}$ and a temporal derivative $\frac{\partial }{\partial t}$:
        \begin{equation}\label{Eq:Transient_PDE}
                    \frac{\partial u({\bf x},t)}{\partial t}  \, = \, {\cal L}(u) \text{ for } {\bf x}\in \Omega \text{ and } 0<t<T_{final}, 
        \end{equation}
    in which $u$ is a scalar function, $\Omega \subset \mathbb{R}^2$ is a closed and bounded spatial domain, and $T_{final}$ is the end time of the simulation.  As our analysis is entirely about the treatment of the unsteady term, the boundary conditions and specific form of ${\cal L}$ are unimportant.  
    
    The TFEM is based on a $N_{el}$ element partition $\Omega^h$ of the domain $\Omega$,  
      \begin{equation}\label{Eq:Partition}
	\Omega^h \, = \, \bigcup_{1\leq k \leq N_{el}} \Omega_k 
      \end{equation}
    in which $\Omega_k$ for all $k$, $1\leq k \leq N_{el}$ are affine transformations, $\left\{\mathcal{J}_k\right\}_{k=1}^{N_{el}}$, of a reference triangular element
    \begin{equation}\label{Eq:Global_Mapping} \mathcal{J}_k(\Omega_{ref}) \, = \, \Omega_k. \end{equation}
    The typical triangular reference element is given by
    \begin{equation}\label{Eq:Reference_Element}
      \Omega_{ref} \, = \, \left\{ (r,s) \, \vert \, -1 \leq r,s \text{ and } r+s \leq 0\right\}.
    \end{equation}
    Locally defined polynomial approximation spaces, denoted $\mathcal{T}_k(p)$ for $1\leq k \leq N_{el}$, are given by mapping 
    \begin{equation}\label{Eq:Triangular_Polynomial_Space}
      \mathcal{T}(p) \, = \, \text{span} \left(\left\{ \left(r^ns^m\right) \, \vert \, 0\leq m,n \text{ and } m+n \leq p\right\}\right)
    \end{equation}
    to the partition elements $\Omega_k$.  The dimension of $\mathcal{T}(p)$ is 
    \begin{equation} \text{dim}\left(\mathcal{T}(p)\right) \, = \, \frac{1}{2}(p+1)(p+2). \end{equation}
            
    In two dimensional PDEs with homogeneous boundary conditions, a continuous finite element method seeks a piecewise polynomial solution restricted to $H^1(\Omega)$, a subset of $C^0(\Omega)$ by the Sobolev Embedding Theorem.
    Let $V_p$ be the finite element space over $\Omega^h$ defined by the local polynomial spaces $ \bigcup\limits_{k=1}^{N_{el}}\mathcal{T}_k(p)$. 
    Define $V^h_p$ as 
    \begin{equation}\label{Eq:Test_and_Trial_Space}
     V^h_p \, = \, V_p \bigcap H^1(\Omega).
    \end{equation}
    then the finite element problem is to find $u^h({\bf x}, t) \in V^h_p$, such that  
       \begin{equation}\label{Eq:Unsteady_Heat_Equation_Variational}
          \int_\Omega \, v^h \frac{\partial u^h({\bf x}, t)}{\partial t} \, d\Omega \, = \, \int_\Omega v^h {\cal L}(u^h({\bf x}, t)) \, d\Omega  \ \ \forall v^h \in V_p^h.
       \end{equation}   
    The space $V^h_p$ is finite dimensional due to the dimensionality of the local approximation spaces, and therefore $u^h$ may be expressed as a linear combination of global basis functions, $\vec{\phi}$, by
       \begin{equation}\label{Eq:Unsteady_Heat_Equation_Spatial_Discretization}
          u^h({\bf x}, t) \, = \, \sum_{}^{} u_i^h(t) \phi_i({\bf x}) \, = \, {\bf U}^T(t) \vec{\phi},
       \end{equation}
    and $v^h$ allows a similar form except with ${\bf U}^T$ replaced by ${\bf V}^T$.
    Eq.~\eqref{Eq:Unsteady_Heat_Equation_Variational} can then be written as
       \begin{equation}\label{Eq:Unsteady_Transient_PDE}
          {\bf V}^T \int_\Omega \vec{\phi} \vec{\phi}^T \, d\Omega \frac{\partial {\bf U}(t)}{\partial t}  \, = \,{\bf V}^T  \int_\Omega \vec{\phi}  {\cal L}\left(\vec{\phi}^T{\bf U}(t)\right) \, d\Omega \,
       \end{equation}
    The global mass matrix, denoted $M$, is defined in the classical way as
       \begin{equation}\label{Eq:Global_Mass}
         M \, = \, \int_\Omega \vec{\phi} \vec{\phi}^T \, d\Omega.
       \end{equation}
    The result is a Galerkin TFEM
       \begin{equation}\label{Eq:TFEM_in_Space}
         M\frac{\partial}{\partial t}{\bf U}(t) \, = \, R({\bf U}(t)).
       \end{equation}
       where $R({\bf U}(t))$ is the discrete form of the spatial operator and the ${\bf V}^T$ term has been removed, as equation~\eqref{Eq:Unsteady_Transient_PDE} must be true for all ${\bf V}^T$.
    
    The model problem \eqref{Eq:Transient_PDE} has thus become a first order differential equation in time:
    \begin{equation}\label{Eq:Explicit_TFEM}
     \frac{d}{d t}{\bf U}(t) \, = \, M^{-1}  R({\bf U}(t)).
    \end{equation}
    In an explicit method the approximate solution at each time step is determined by a matrix multiplication of the inverse mass matrix $M^{-1}$.
    To be computationally efficient, it is therefore necessary to have a mass matrix that is trivially inverted, i.e. a diagonal mass matrix.
    
    \begin{subsection}{High Order Bases for Continuous Methods}\label{Sec:Continuity}
    Following along with the work of Dubiner in \cite{Dubiner1991}, we present necessary conditions for the direct enforcement of a globally $C^0$ approximation space.
      Specifically, the reference element $\Omega_{ref}$ is further decomposed into the vertices $V_A$, $V_B$, $V_C$, and the edges $E_{\mathcal A}$, $E_{\mathcal B}$, and $E_{\mathcal C}$ as shown in Fig.~\ref{Fig:Spaces}.
      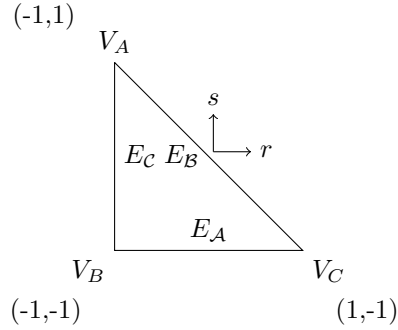
\begin{figure}[H]
	\centering
	  \def\axislength{.2}
        \begin{tikzpicture}[scale=1.25]

        \coordinate (0) at (0.05,0.05);
        \coordinate (r) at (0.45,0.05);
        \coordinate (s) at (0.05,0.45);
        \coordinate [label=above:{(-1,1)}] (lt) at (-1.75, 1.25);
        \coordinate  (lb) at (-1.75, -1.15);
        \coordinate [label=left:{(-1,-1)}] (bl) at (-1.25, -1.65);
        \coordinate [label=right:{(1,-1)}] (br) at (1.25, -1.65);
        \coordinate [label=above:$V_A$] (v1) at (-1,1);
        \coordinate [label=below left:$V_B$] (v2) at (-1,-1);
        \coordinate [label=below right:$V_C$] (v3) at (1,-1);
        \draw (v1) -- node[right] {$E_{\cal C}$} (v2) -- node[above] {$E_{\cal A}$} (v3) -- node[left] {$E_{\cal B}$} (v1);
        \draw[->] (0) -- (r) node[right] {$r$};
        \draw[->] (0) -- (s) node[above] {$s$};

        \end{tikzpicture}
	  \caption{The vertex and edge arrangement of the reference triangle $\Omega_{ref}$.}\label{Fig:Spaces}
      \end{figure}
      \noindent The concept of vertex, edge, and interior spaces is used to categorize local basis functions:
	A \textit{vertex function} is a function $f:\Omega_{ref}\rightarrow\mathbb{R}$ with non-zero evaluation, typically $1$, at only one of the vertices and zero evaluation along the opposing edge.
	An \textit{edge function} is a function $f:\Omega_{ref} \rightarrow \mathbb{R}$ with non-zero evaluation along only one edge and zero evaluation along the other two edges.
	And, an \textit{interior function} is a non-zero function $f:\Omega_{ref}\rightarrow \mathbb{R}$ such that the evaluation of $f$ is zero along all edges (and at all vertices).
      For a basis for the space ${\cal T}(p)$ to be useful for $C^0$ finite elements, there must be one vertex function for each vertex, $(p-1)$ edge functions for each edge, and $(p-1)(p-2)/2$ interior functions.
      
      In \cite{Dubiner1991}, Dubiner introduced a basis, $\vec{\phi}$ for $\mathcal{T}(p)$ that is commonly used for $C^0$ TFEM.
      This basis is achieved by mapping $\Omega_{ref}$ to the reference quadrilateral element, $\{(\xi, \eta)\vert-1\leq \xi, \eta \leq 1\}$, by
    \begin{equation}\label{Eq:Dubiner_Mapping}
      \xi = -1+2\left( \frac{1+r}{1-s}\right) \text{ and } \eta = s.
    \end{equation}
    The reference triangle can then be represented in terms of a warped coordinate system $(\xi, \eta)$, see Fig.~\ref{Fig:Tris}.
    \begin{figure}[H]
      \centering
      \includegraphics[width=.5\textwidth]{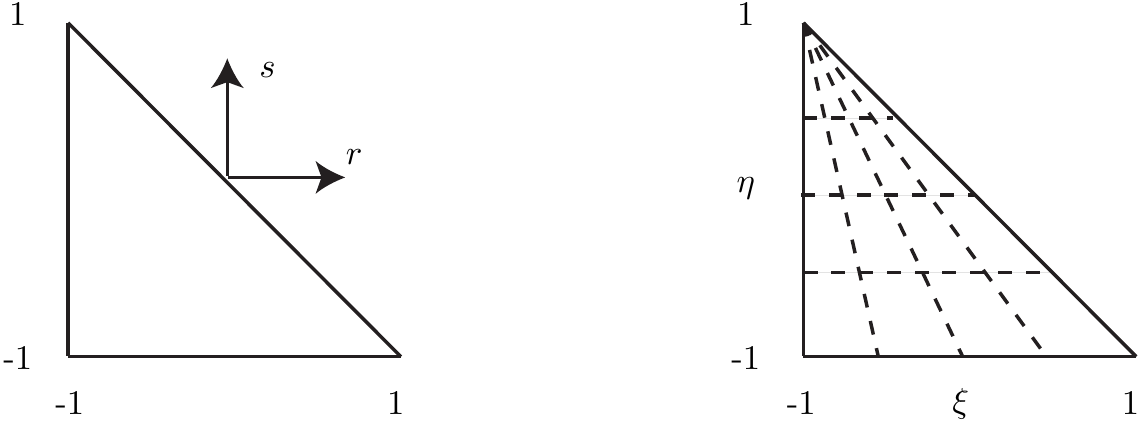}
      \caption{The reference triangle $\Omega_{ref}$ shown in terms of $(\xi, \eta)$.}\label{Fig:Tris}
    \end{figure}
    
    In terms of the warped coordinates $(\xi, \eta)$, the modified Dubiner basis functions are originally defined in \cite{Dubiner1991}, then slightly modified in \cite{Helenbrook2009}.
    The three vertex functions are given by 
      \begin{equation}\label{Eq:Dubiner_Vertex_Functions} 
	\begin{array}{c} 
	  \phi_{A}(\xi, \eta) \, = \, \left(\frac{1+\eta}{2}\right), \\  
	  \phi_{B}(\xi, \eta) \, = \, \left(\frac{1-\xi}{2}\right)\left(\frac{1-\eta}{2}\right), \\ 
	  \phi_{C}(\xi, \eta) \, = \, \left(\frac{1+\xi}{2}\right)\left(\frac{1-\eta}{2}\right).
	\end{array}
      \end{equation}
      For $p\geq 2$ the basis function definitions rely on the classical Jacobi Polynomials $P^{\alpha, \beta}_n(z)$, these polynomials are orthogonal over the interval $(-1,1)$ with respect to the weighting function $(1-z)^\alpha(1+z)^\beta$. 
      The $(p-1)$ edge functions for each edge are
      \begin{equation}\label{Eq:Dubiner_Edge_Functions} 
	\begin{array}{c c}
	  \phi_{{\cal A}, m} \, = \,  \left( \frac{1+\xi}{2} \right) \left( \frac{1-\xi}{2} \right) P_{m-1}^{2,2}(\xi) \left( \frac{1-\eta}{2} \right)^{m+1}    & \text{ for } 1\leq m \leq p-1,  \\
	  \phi_{{\cal B}, m} \, = \,  \left( \frac{1+\xi}{2} \right) \left( \frac{1-\eta}{2} \right)  \left( \frac{1+\eta}{2} \right) P_{m-1}^{2,2}(\eta)       & \text{ for } 1\leq m \leq p-1, \\
	  \phi_{{\cal C}, m} \, = \,  (-1)^{m-1}\left( \frac{1-\xi}{2} \right) \left( \frac{1-\eta}{2} \right)  \left( \frac{1+\eta}{2} \right) P_{m-1}^{2,2}(\eta) & \text{ for } 1\leq m \leq p-1.
	\end{array}
      \end{equation}
      The interior functions for $0\leq m \leq (p-3)$ and $0 \leq n \leq (p-3-m)$ are defined
      \begin{equation}\label{Eq:Dubiner_Interior_Functions} 
	  \phi_{{\cal I}, m, n} \, = \, \left( \frac{1-\xi}{2} \right) \left( \frac{1+\xi}{2} \right) P_m^{2,2}(\xi) \left( \frac{1-\eta}{2} \right)^{m+2}  \left( \frac{1+\eta}{2} \right) P_n^{2m+5,2}(\eta) \\
      \end{equation}
      and given one-dimensional indices using
      \begin{equation}\label{Eq:Dubiner_Interior_Ordering}
	j(m,n) \, = \, \frac{1}{2}(m+n)(m+n+1)+n+1.  
      \end{equation}
      
      An important property of this basis, and linearly independent sets of polynomials in general, regarding subsets and their linear independence is given in the following lemma which states specifically that the linear independence of a set of polynomials is not affected by the removal of common factors.   
      This property will later be used to show the invertibility of a matrix with elements defined by the $L_2$ inner product.
        
      \begin{lemma}\label{Lemma:Basis_Functions_mod_f}
        \textit{Any finite dimensional, linearly independent set of polynomials possessing a common factor $f$, $\{f g_i\}$ defines the linearly independent set, $\{g_i\}$.} 
          \begin{proof}
	    To prove this consider the contrapositive:
	    \textit{The linear dependence of a set of polynomials is maintained after element-wise multiplication of any function. }
            Let $G$ be a set $N$ of linearly dependent functions over some space $\Omega$,
                \[ G \, := \, \left\{g_1, g_2, \dots, g_N \right\}. \]
            Then there exists some $n$, $1\leq n \leq N$, for which $g_n \in G$ is a linear combination of the other functions:
                \[ g_n \, = \, \sum_{k=1, \, k\neq n}^N a_k g_k\]
            Let $f$ be any function defined on $\Omega$. 
            Define the set of functions $\hat{G}$,
                \[ \hat{G} \, := \, \left\{ fg_1, fg_2, \dots, fg_N\right\}. \]
            It is then clear that
                \[ fg_n \, = \, f\left( \sum_{k=1, \, k\neq n}^N a_k g_k \right)\, = \,\sum_{k=1, \, k\neq n}^N a_k f g_k. \]    
            Therefore, $\hat{G}$ is linearly dependent as well. 
        \end{proof}   
    \end{lemma}
    As an example application of this lemma, the set of edge $E_{{\cal B}}$ basis polynomials for $\mathcal{T}(p)$ all possess common factors of $\phi_{A}$ and $\phi_{C}$, the vertex $V_A$ and $V_C$ functions.
    The factor $\phi_{A}$ ensures that the functions $\phi_{{\cal B}, i}$ for $1 \leq i \leq p-1$ are all zero along $E_{{\cal A}}$, and the factor $\phi_{C}$ forces a zero evaluation along $E_{\cal C}$. 
    The functions $\phi_{{\cal B}, i}$ are polynomials of degree varying between $2$ and $p$, whereas $\phi_{A}$ and $\phi_B$ are first degree polynomials.
    By factoring $\phi_{A}$ and $\phi_{C}$ out of the edge functions, $\phi_{\cal B} / \left( \phi_A \phi_{C} \right)$, the resulting set of polynomials remains linearly independent with degree ranging between $0$ and $p-2$.   
    In fact, the resulting set is a linearly independent subset of $\mathcal{T}(p-2)$.
    
    \end{subsection}   
 \end{section}


   \begin{section}{Diagonal pseudo-mass matrix}\label{Sec:Diagonal_Pseudo_Mass}
      The $L_2$ projection of a function $f$ in $L_2(\Omega_{ref})$ onto $\mathcal{T}(p)$ is the unique function in $\mathcal{T}(p)$ that minimizes the error measured in the $L_2$ norm.
      In finite element methods the $L_2$ projection is evaluated by first identifying the coefficient vector of the projection and then taking the inner product of this coefficient vector with the basis.
      For any basis $\vec{\psi}$ of $\mathcal{T}(p)$ over $\Omega_{ref}$ there exists an invertible mass matrix $M_{\psi}$ defined by \eqref{Eq:Global_Mass}.
      The coefficient vector of the $L_2$ projection is 
      \begin{equation}\label{Eq:Coefficient_Vector}
	 \vec{u}_{\psi} \, = \, M_{\psi}^{-1} \,\int_{\Omega_{ref}}  \, f \vec{\psi} \, d\Omega. 
      \end{equation}
      The unique $L_2$ projection is then given by $\vec{u}_{\psi}^T \vec{\psi}$, and for all $f \in \mathcal{T}(p)$ the $L_2$ projection of $f$ is identically $f$:
      \begin{equation}\label{Eq:Unique_Rep}
	\vec{u}_{\psi}^T \vec{\psi} \, = \, f \text{ for all } f \in \text{span}\left(\vec{\psi}\right).
      \end{equation} 

    There is no basis that is both orthogonal and divisible into vertex, edge, and interior functions. This implies that there is no basis for $C^0$ FEM with a diagonal mass matrix.   
    There are  however \textit{pseudo-mass matrix} approaches; a pseudo-mass matrix is simply a non-singular real-valued matrix $\underbar{M} \in \mathbb{R}^{\text{dim}(\mathcal{T}(p))\times\text{dim}(\mathcal{T}(p))}$ that acts in place of the mass matrix.
    Mass-lumping techniques, which replace the mass matrix with a row-summed diagonal matrix, clearly constitute pseudo-mass matrix approaches \cite{Armentano2003, Cohen2001, Durufle2009, Franca1997, Kong1999}.  However, these approaches often suffer from poor accuracy.
    
    SEM is a pseudo-mass approach as the approximate mass matrix is cleverly generated by the GL integration scheme at the expense of only one polynomial degree of accuracy.  In hopes of matching this success, we seek a high-order method that uses an accurate diagonal pseudo-mass matrix.
    The $p^{th}$ degree SEM pseudo-mass matrix is capable of exactly projecting all polynomials of one less polynomial degree than the basis, namely polynomials of degree $(p-1)$.
    We describe the accuracy of the pseudo-mass matrix by the degree of the subspace for which the projections are exact.
    For a FEM using a $p$ degree basis, the pseudo-mass matrix $\underbar{M}$ is said to be $k$-exact for $1\leq k \leq p$ when the coefficient vector
    \begin{equation*}
    \vec{u}_\psi \, = \, \underbar{M}^{-1} \left( \int_\Omega \, \vec{\psi} f \, d\Omega\right)
    \end{equation*}
    is exact for all $f$ of polynomial degree up to $k$;
    for TFEM
    \[\vec{u}_\psi^T \vec{\psi} \, = \, f \text{ for all } f\in \mathcal{T}(k) \subseteq \mathcal{T}(p).\]
    The SEM uses a $(p-1)$-exact pseudo-mass matrix and it would be ideal for the TFEM to have a basis $\vec{\psi}$ for which there exists a $(p-1)$-exact diagonal pseudo-mass matrix, $D$: 
    \begin{equation}\label{Eq:D_Pseudo_Mass}
    D \, = \, \left[
    \begin{array}{cccc}
    D_1 & & & \\
    & D_2 & & \\
    & & \ddots & \\
    & & & D_{\text{dim}\left(\mathcal{T}(p)\right)}   
    \end{array}
    \right].
    \end{equation}     
    
    \begin{subsection}{An Arbitrary Change of Basis for Continuous Methods}\label{Sec:Change_of_Basis}
        If a $(p-1)$-exact diagonal pseudo-mass matrix exists, the basis for which it functions is unknown.
        Therefore, it is beneficial to identify all sets of bases appropriate for the $C^0$ TFEM.
        Recall the modified Dubiner basis $\vec{\phi}$, \eqref{Eq:Dubiner_Vertex_Functions}~--~\eqref{Eq:Dubiner_Interior_Functions}.
        As $\vec{\phi}$ is a basis for $\mathcal{T}(p)$, for any function $f\in \mathcal{T}(p)$,  by \eqref{Eq:Unique_Rep} $\vec{u}_\phi$ defined via \eqref{Eq:Coefficient_Vector} satisfies
        \[ \vec{u}_{\phi}^T \vec{\phi} \, = \, f.\]
        Furthermore by the equivalence of bases, any basis $\vec{\psi}$ may be written as a change of basis $T \in \mathbb{R}^{\text{dim}(\mathcal{T}(p))\times\text{dim}(\mathcal{T}(p))}$ of the modified Dubiner basis $\vec{\phi}$:      
        \begin{equation}\label{Eq:New_Basis}
            \vec{\psi} \, = \, T\vec{\phi}.
        \end{equation}
        
        Any basis $\vec{\psi}$ need only be defined by an invertible $T$ coupled with $\vec{\phi}$.
        However, an arbitrary invertible $T$ does not guarantee that the basis $\vec{\psi}$ satisfies the continuity constraints. 
        Fortunately, the modified Dubiner basis is a continuous basis.
        Since the modified Dubiner basis satisfies the continuity constraints, the mapping $T$ must only enforce the preservation of these constraints in the new basis.
        The continuity constraints given earlier, \S\ref{Sec:Continuity}, state the need of vertex, edge, and interior basis functions.

        Let the vertex $V_A$, vertex $V_B$, and vertex $V_C$ basis functions be represented by the subscripts $A$, $B$, and $C$ respectively. 
        For example the vertex $V_A$ function of $\vec{\phi}$ is written $\phi_{A}$.
        There will be $(p-1)$ functions along each edge, and these will be denoted by the subscripts ${\cal A}_i$, ${\cal B}_i$, and ${\cal C}_i$ for $1\leq i \leq (p-1)$. 
        The functions $\phi_{{\cal A}_i}$, $\phi_{{\cal B}_i}$, and $\phi_{{\cal C}_i}$ will be non-zero on the edges $E_{\cal A}$, $E_{\cal B}$, and $E_{\cal C}$ respectively (shown in Fig.~\ref{Fig:Spaces}).
        Finally, ${\cal I}_i$ will be used to specify the interior functions, organized by $i = j(m,n)$, \eqref{Eq:Dubiner_Interior_Ordering}.
        
        The vertex functions of $\vec{\phi}$ must remain vertex functions under transformation by $T$;  the first three rows of $T$ define the three vertex functions of the new basis.
        Recall that any vertex $V_A$ function, $\psi_A$ must be zero along the edge $E_{\cal A}$.
        Therefore $\psi_A$ exists in the subspace of $\mathcal{T}(p)$ spanned by basis functions of $\vec{\phi}$ that are zero on $E_{\cal A}$: $\phi_A$, $\phi_{{\cal B}, i}$, $\phi_{{\cal C}, i}$, and $\phi_{{\cal I}, i}$.
        Any new vertex $V_A$ function $\psi_A$ is then a linear combination of these functions:
        \[\psi_{A} \, = \, \phi_A + \sum_{i=1}^{p-1} A_{{\cal B}, i}\phi_{{\cal B}, i} + \sum_{i=1}^{p-1}A_{{\cal C},i}\phi_{{\cal C}, i} + \sum_{j=1}^{\frac{1}{2}(p-1)(p-2)}A_{{\cal I}, j} \phi_{{\cal I}, j},\]
        in which $\psi_A$ has been normalized by $\phi_A$ and the coefficients $A_{{\cal B}, i}$, $A_{{\cal C},i}$, and $A_{{\cal I}, j}$ represent the additions of edge $E_{\cal B}$, edge $E_{\cal C}$, and interior functions to $\phi_A$ respectively.
        The first three rows of $T$ are then given by
        \begin{equation}\label{Eq:T_First_Three_Rows}
        T_{1:3, :} \, = \, \left[\begin{array}{ccccccc}
                                    1 & & & & A_{\cal B} & A_{\cal C} & A_{\cal I} \\
                                    & 1 & & B_{\cal A} & & B_{\cal C} & B_{\cal I} \\
                                    & & 1 & C_{\cal A} & C_{\cal B} & & C_{\cal I}
                                \end{array} \right],
        \end{equation}
        for which $A_{\cal B}$, $A_{\cal C}$, $B_{\cal A}$, $B_{\cal C}$, $C_{\cal A}$, and $C_{\cal B}$ are all scalar valued row vectors of length $(p-1)$, and the interior augmentations are row vectors of length $(p-1)(p-2)/2$.
        
        The next $3(p-1)$ rows of $T$ are used to map the edge functions of $\vec{\phi}$.
        The modified Dubiner basis edge functions are organized first by edge and then by degree such that basis functions $4$ to $4+(p-1)$ are edge $E_{\cal A}$ functions, followed by the edge $E_{\cal B}$ functions, then $E_{\cal C}$ functions.   
        Note that this organization of basis functions is used in this section to place emphasis on the structure of individual edge subspaces of $\mathcal{T}(p)$, it will be altered in a later section.
        By definition, an edge function exists in a space spanned by functions that are non-zero along only one edge.
        The same-edge edge functions and the interior functions span this subspace. 
        Similarly, only the interior functions satisfy the interior function requirements.
        Therefore, the most general transformation from the modified Dubiner basis to any other basis adhering to the continuity constraints is given by 
        \begin{equation}\label{Eq:Transformation_Shape}
                T \, = \, \left[
                    \begin{array}{c c c c c c c}
                    1	&	&	& 			& A_{\cal B} & A_{\cal C}		&A_{\cal I} \\
                    & 1	&	& B_{\cal A} 	&			& B_{\cal C}	&B_{\cal I} \\
                    &	& 1	& C_{\cal A} & C_{\cal B}&				& C_{\cal I}  \\
                    &  &  & {\cal A}_{\cal{A}} &  &  & {\cal A}_{\cal I} \\
                    &  &  &  & {\cal B}_{\cal{B}} &  &  {\cal B}_{\cal I} \\
                    &  &  &  &  & {\cal C}_{\cal{C}} &  {\cal C}_{\cal I} \\
                    &  &  &  &  &  & {\cal I}  \\
                    \end{array}
                \right],
            \end{equation}
        in which ${\cal A}_{\cal A}$, ${\cal B}_{\cal B}$, ${\cal C}_{\cal C}$, and ${\cal A}_{\cal I}$, ${\cal B}_{\cal I}$, and ${\cal C}_{\cal I}$ denote the additions of same-edge and interior functions.
        The blocks ${\cal A}_{\cal A}$, ${\cal B}_{\cal B}$, and ${\cal C}_{\cal C}$ are all $(p-1)\times (p-1)$ non-singular matrices.   
        The dimensions of ${\cal A}_{\cal I}$, ${\cal B}_{\cal I}$, and ${\cal C}_{\cal I}$ are all $(p-1)\times (p-1)(p-2)/2$; and the $(p-1)(p-2)/2\times(p-1)(p-2)/2$  non-singular block ${\cal I}$ is used for the interior change of basis.
        As all diagonal blocks of $T$ are non-singular, $T$ is invertible.
        Thus, any basis $\vec{\psi}$ satisfying the continuity constraints can be represented by $\vec{\phi}$ via $T \vec{\phi}$ with $T$ restricted to the structure of \eqref{Eq:Transformation_Shape}.
    \end{subsection}
    
    \begin{subsection}{No \texorpdfstring{$(p-1)$}--exact Diagonal TFEM}\label{Sec:Proof}   
   From \cite{Helenbrook2009}, there are unique vertex functions for $p\geq1$ that can be coupled with a $(p-1)$-exact diagonal pseudo-mass matrix.
   Therefore the first three rows of the transformation matrix exist.  It remains to be determined if edge functions and interior functions can be found.
   Assume that the entire $(p-1)$-exact pseudo-mass matrix $D$ and basis $\vec{\psi}$ exist.
      So, for any function $f \in \mathcal{T}(p-1)$ there exists a unique $\vec{u}_\psi$ defined by 
      \begin{equation}\label{Eq:Diagonal_System} \vec{u}_\psi \, = \, D^{-1} \int f \vec{\psi} \, d\Omega \end{equation}
      such that $\vec{u}_{\psi}^T\vec{\psi} \, = \, f$.
      This differs from \eqref{Eq:Coefficient_Vector} because $D$ is not the mass matrix, and it is not exact for functions in $\mathcal{T}(p)/\mathcal{T}(p-1)$.
      
   As the three vertex functions have been identified by Helenbrook \cite{Helenbrook2009}, we look more closely to the first edge function $\psi_{{\cal A},1}$.
	Associated with the fourth row of \eqref{Eq:Diagonal_System}, this edge function is defined in \eqref{Eq:New_Basis} by the following linear combination of modified Dubiner basis functions:
	\begin{equation}\label{Eq:First_First_Edge_Mode}
	  \psi_{{\cal A}, 1} \, = \, \sum_{j=1}^{p-1} {\cal A}_{{\cal A},1,j} \phi_{{\cal A},j} + \sum_{k=1}^{(p-1)(p-2)/2} {\cal A}_{{\cal I},1,k} \phi_{{\cal I},k}.
	\end{equation}
	By assumption, the fourth row of \eqref{Eq:Diagonal_System} should hold for all functions $f\in \mathcal{T}(p-1)$.   These constraints will be identified by means of the modified Dubiner basis. 
 
    Let $f \in \mathcal{T}(p-1)$ and let $\vec{\psi}$ be a basis for which a $(p-1)$-exact diagonal pseudo-mass matrix exists.
      Note that because $f \in \mathcal{T}(p-1) \subset \mathcal{T}(p)$, there is a unique coefficient vector defined in \eqref{Eq:Coefficient_Vector} for which
          \[  \vec{u}_\phi^T\vec{\phi} \, = \,  \vec{u}_\psi^T \vec{\psi} \, = \,  f. \]
      From \eqref{Eq:New_Basis}, $\vec{\psi} = T\vec{\phi}$, and $ \vec{u}_\psi^T \left(T\vec{\phi}\right) \, = \, \vec{u}_\phi^T\vec{\phi},$ or
        \[ \left(\vec{u}_\psi^T T\right)\vec{\phi} \, = \, \vec{u}_\phi^T\vec{\phi}.\]
      As the coefficient vectors are unique, $ T^{T}\vec{u}_\psi \, = \, \vec{u}_\phi.$
      The coefficient vector for $\vec{\psi}$ is then
                \begin{equation}\label{Eq:Psi_Soln_by_Phi}
           \vec{u}_\psi \, = \,  T^{-T} \vec{u}_\phi.
          \end{equation}
      The inverse of $T$ can be found by block inversion of \eqref{Eq:Transformation_Shape}:
      \begin{equation}\label{Eq:Inverse_Transformation_Shape}
                T^{-1}  =  \left[
                    \begin{array}{c c c c c c c}
                    	1	&	&	& 			& -A_{\cal B} {\cal B}_{\cal B}^{-1} & -A_{\cal C} {\cal C}_{\cal C}^{-1}	& \left[ T^{-1} \right]_{A,{\cal I}} \\
                    	& 1	&	& -B_{\cal A} {\cal A}_{\cal A}^{-1}&			& -B_{\cal C} {\cal C}_{\cal C}^{-1}                    & \left[ T^{-1} \right]_{B,{\cal I}} \\
                    	&	& 1	& -C_{\cal A} {\cal A}_{\cal A}^{-1} & -C_{\cal B} {\cal B}_{\cal B}^{-1} &				        & \left[ T^{-1} \right]_{C,{\cal I}} \\
                    	&  &  & {\cal A}_{\cal A}^{-1} &  &  & -{\cal A}_{\cal A}^{-1}{\cal A}_{\cal I}\\
                    	&  &  &  & {\cal B}_{\cal B}^{-1} &  & -{\cal B}_{\cal B}^{-1}{\cal B}_{\cal I} \\
                    	&  &  &  &  & {\cal C}_{\cal C}^{-1} & -{\cal C}_{\cal C}^{-1}{\cal C}_{\cal I} \\
                    	&  &  &  &  &  & {\cal I}^{-1} \\
                    \end{array}\right]
            \end{equation}   
        with 
        \[ \left[ T^{-1} \right]_{A,{\cal I}} \, = \, -A_{\cal I}+A_{\cal B} {\cal B}_{\cal B}^{-1} {\cal B}_{\cal I} + A_{\cal C} {\cal C}_{\cal C}^{-1} {\cal C}_{\cal I},\]
        \[ \left[ T^{-1} \right]_{B,{\cal I}} \, = \, -B_{\cal I}+B_{\cal A} {\cal A}_{\cal A}^{-1} {\cal A}_{\cal I} + B_{\cal C} {\cal C}_{\cal C}^{-1} {\cal C}_{\cal I},\]
        and 
        \[ \left[ T^{-1} \right]_{C,{\cal I}} \, = \, -C_{\cal I} +C_{\cal A} {\cal A}_{\cal A}^{-1} {\cal A}_{\cal I} +C_{\cal B} {\cal B}_{\cal B}^{-1} {\cal B}_{\cal I}.\]
            
	The constraints will consist of the $L_2$ projections of the basis functions of $\mathcal{T}(p-1)$. As the modified Dubiner basis is hierarchical, $f$ in \eqref{Eq:Diagonal_System} can be replaced by $\phi_i$ for all indices $i$ such that $\phi_i\in\mathcal{T}(p-1)$.  This in turn implies that $\vec{u}_\phi$ is zero at every entry except for the $i^{th}$ entry which is one.  Substituting \eqref{Eq:First_First_Edge_Mode} and \eqref{Eq:Psi_Soln_by_Phi} on the right and left hand side of \eqref{Eq:Diagonal_System} respectively, the fourth row of \eqref{Eq:Diagonal_System} can be written as 
	\begin{equation}\label{Eq:Contradiction_Constraints}
	 D_{4,4}T^{-T}_{4,i} \, = \, \int_\Omega \, \phi_i\left( \sum_{j=1}^{p-1} {\cal A}_{{\cal A},1,j} \phi_{{\cal A},j} + \sum_{k=1}^{(p-1)(p-2)/2} {\cal A}_{{\cal I}, 1, k} \phi_{{\cal I}, k}\right) \, d\Omega.
	\end{equation}
	
	Equation \eqref{Eq:Contradiction_Constraints} is a system of equations and may be more easily solved by isolating the unknowns from ${\cal A}_{\cal A}$ and ${\cal A}_{\cal I}$ on the right hand side of the equation.
	Denote these unknowns as the vector $\vec{T}^*_4$ defined as
	\begin{equation}\label{Eq:Vec_T_Star}
	 \vec{T}^*_4 \, = \, \left[ \begin{array}{c}  
		    {\cal A}_{{\cal A},1,1} \\ 
		    {\cal A}_{{\cal A},1,2} \\ 
		    \vdots \\
		    {\cal A}_{{\cal A},1,p-1} \\
		    {\cal A}_{{\cal I},1,1} \\
		    {\cal A}_{{\cal I},1,2} \\ 
		    \vdots \\
		    {\cal A}_{{\cal I},1,(p-1)(p-2)/2} \\
	         \end{array} \right].
	\end{equation}
	
	Using \eqref{Eq:Vec_T_Star}, \eqref{Eq:Contradiction_Constraints} is
	\begin{equation}\label{Eq:Contradiction_Constraints_3}
	  D_{4,4}T^{-T}_{4,i} \, = \, \int_\Omega \, \phi_i\left[ \phi_{{\cal A},1}:\phi_{{\cal A},p-1}, \, \phi_{{\cal I}, 1}:\phi_{{\cal I}, (p-1)(p-2)/2}\right] d\Omega \, \vec{T}^*_4
	\end{equation}
	for all $i$ such that $\phi_i \in  \mathcal{T}(p-1)$.  
	In \S\ref{Sec:Proof_indices} we formally state, and reduce, a list of indices $i$ of constraining functions.
	Then in \S\ref{Sec:Proof_Proof} we prove that this system results in $\vec{T}^*_4$ being identically zero which contradicts the fact that $T$ must be invertible.

	\begin{subsubsection}{Defining the Constraining Indices}\label{Sec:Proof_indices}
	The constraints are all of the functions of one less polynomial degree, $(p-1)$.
	Recall that the indices are currently arranged as vertices, $[A,B,C]$, edge ${\cal A}$, ${\cal B}$, and ${\cal C}$; then finally the interior functions, ${\cal I}$.
	The edge and interior functions are organized by increasing order.
	For $p \geq 2$, which is the case of interest here, vertex functions are always part of the constraints, which implies that the indices $[1,2,3]$, representing $[A,B,C]$, are included.   
	For $p > 2$ all edge functions of $\vec{\phi}$ exist in $\mathcal{T}(p-1)$ except for the last function on each side.  
	Thus the indices $\left[4\!:\!(p+1), (p+3)\!:\!2p, (2p+2)\!:\!(3p-1)\right]$ are included.  
	Lastly, for $p>3$ all interior functions of  $\mathcal{T}(p-1)$  are also a component of $\vec{\phi}$.  
	The indexing system, \eqref{Eq:Dubiner_Interior_Ordering}, for the interior functions was created such that as $p$ increases, functions are added to the end of the indexing system, thus the interior constraints are the first $(p-2)(p-3)/2$ interior basis functions.  
	This results in the following list of indices
	\begin{equation}\label{Eq:mu}
	\mu:=  \left[1\!:\!3, 4\!:\!(p+1), (p+3)\!:\!2p, (2p+2)\!:\!(3p-1), (3p+1)\!:\!(3p+(p-2)(p-3)/2) \right].
	\end{equation}
	The inverse of a matrix transposed is the transpose of the inverse of that matrix; rewriting $T^{-T}_{4,i}$ as $T^{-1}_{i,4}$ in \eqref{Eq:Contradiction_Constraints_3} gives
	\begin{equation}\label{Eq:Contradiction_Constraints_2}
	  D_{4,4}T^{-1}_{i,4} \, = \, \int_\Omega \, \phi_i\left[ \phi_{{\cal A},1}:\phi_{{\cal A},p-1}, \, \phi_{{\cal I}, 1}:\phi_{{\cal I}, (p-1)(p-2)/2}\right] d\Omega \, \vec{T}^*_4\ \forall i \in \mu.
	\end{equation}

	  The desired basis $\vec{\psi}$ is unknown, and therefore $T$ is unknown as well.  However, a complete knowledge of the fourth column of $T^{-1}$ is unnecessary in the derivation of a contradiction. From  \eqref{Eq:Inverse_Transformation_Shape},  $D_{4,4}T^{-1}_{i,4}$ has only $(p+1)$ non-zero entries, namely one from the term $-B_{\cal A} {\cal A}_{\cal A}^{-1}$, denoted as $T^{-1}_{2,4}$, one from the term $-C_{\cal A} {\cal A}_{\cal A}^{-1} $ denoted as  $T^{-1}_{3,4}$, and $p-1$ from the term ${\cal A}_{\cal A}^{-1}$ denoted as $T^{-1}_{4:(p+2), 4}$.	The rows of $T^{-1}_{i,4}$ associated with the first vertex function, the edge two and three functions, and the interior functions are zero. 
	The subset of indices that result in a zero value for $D_{4,4}T^{-1}_{i,4}$ will be referred to as the zero constraints of the fourth row for $p$, denoted as $\mu_0$:  
	\begin{equation}\label{Eq:Zero_Constraints}
	 \mu_0 := \left[1, (p+3)\!:\!2p, (2p+2)\!:\!(3p-1), (3p+1)\!:\!(3p+(p-2)(p-3)/2) \right].
	\end{equation}
	Then for $i \in \mu_0$,
	\begin{equation} \label{Eq:Contradiction_Constraints_4}
	\, \int_\Omega \, \phi_i\left[ \phi_{{\cal A},1}:\phi_{{\cal A},p-1}, \, \phi_{{\cal I}, 1}:\phi_{{\cal I}, (p-1)(p-2)/2}\right] d\Omega \, \vec{T}^*_4 = 0.
	\end{equation}
	\end{subsubsection}
	
	\begin{subsubsection}{No Accurate Diagonal Scheme Exists}\label{Sec:Proof_Proof}
	
	If there exist sufficient non-degenerate constraints, an invertible linear system defining the elements of the fourth row of $T$ and $D_{4,4}$ is defined.	Adding up the number of entries in $\mu_0$ determines that
	    \[ \text{dim}\left(\mu_0\right) \, = \, 1+2(p-2)+(p-2)(p-3)/2 \, = \, \frac{1}{2}p(p-1).\]
	Similarly, $\vec{T}^*_4$ as defined by \eqref{Eq:Vec_T_Star} has $(p-1)$ degrees of freedom from the first row of ${\cal A}_{\cal A}$ and $(p-1)(p-2)/2$ degrees of freedom from ${\cal A}_{\cal I}$ for a total of
	\[    \text{dim}\left(\vec{T}^*_4\right) =  (p-1)+\frac{1}{2}(p-1)(p-2) \, = \, \frac{1}{2}p(p-1) \]
	degrees of freedom.
	
	\begin{theorem}{The system defined by \eqref{Eq:Contradiction_Constraints_4} for $i \in  \mu_0$ is both square and non-singular. }\label{Theorem:Fourth_Row_of_T_System}
                
                \begin{proof}
                     Let $\vec{\phi}_{\mu_0}$ denote the set of basis functions from $\vec{\phi}$ that force $D_{4,4}T^{-1}_{i,4} = 0$,  i.e. $\phi_i$ for $i \in \mu_0$. 
                    Let $\vec{\phi}_{*}$ denote the set of basis functions of $\vec{\phi}$ associated with $\vec{T}^*_4$ as shown in \eqref{Eq:Contradiction_Constraints_4}.  As discussed above,  both $\vec{\phi}_{\mu_0}$ and $\vec{\phi}_{*}$ have $\frac{1}{2}p(p-1)$ elements. 
                    The elements of $\vec{\phi}_{\mu_0}$ share a common factor $\phi_{A}$ because all of the functions are zero along edge ${\cal A}$.  This can be seen by referring to \eqref{Eq:Dubiner_Vertex_Functions}, \eqref{Eq:Dubiner_Edge_Functions}, and \eqref{Eq:Dubiner_Interior_Functions}.
                    Similarly, every element of $\vec{\phi}_*$ shares the common factors $\phi_{B}$ and $\phi_{C}$ because the function is required to be zero along edges ${\cal B}$ and ${\cal C}$.  This can also be verified by direct examination of the basis.    
                    Define $\hat{\phi}_{\mu_0}$ and $\hat{\phi}_*$ such that
                        \begin{equation}\label{Eq:Reduced_alpha}
			  \vec{\phi}_{\mu_0} \, = \, \phi_{A} \hat{\phi}_{\mu_0},
			\end{equation}
                    and   
                        \begin{equation}\label{Eq:Reduced_beta}
			  \vec{\phi}_* \, = \, \phi_{B}\phi_{C} \hat{\phi}_*.
			\end{equation}
                    
                    The elements of $\vec{\phi}_{\mu_0}$ range in degree from $1$ through $p-1$.
                    A polynomial of degree $1$ was factored from this.  
                    This leaves us with $\frac{1}{2}p(p-1)$ polynomials of degree $0$ through $p-2$, or consequently $\text{dim}\left(\mathcal{T}(p-2)\right)$.
                    By Corollary~\ref{Lemma:Basis_Functions_mod_f}, these are linearly independent in $\mathcal{T}(p-2)$.
                    So, the elements of $\hat{\phi}_{\mu_0}$ form a basis for $\mathcal{T}(p-2)$.                  
                    Similarly, the elements of $\vec{\phi}_*$ range in degree from $2$ through $p$.
                    Out of $\vec{\phi}_*$ two vertex functions were factored, leaving the elements of $\hat{\phi}_*$ ranging in degree $0$ through $p-2$.
                    Therefore, the elements of $\hat{\phi}_*$ form a basis for $\mathcal{T}(p-2)$ as well.
                    
                    By the equivalence of basis functions, there exists an invertible transformation, not necessarily of similar shape to \eqref{Eq:Transformation_Shape}, but invertible none-the-less, that maps between $\hat{\phi}_{\mu_0}$ and $\hat{\phi}_*$.
                    Call the transformation $H$.
                    Therefore
                        \[ \hat{\phi}_{\mu_0} \, = \, H \hat{\phi}_*.\]
                    
                    The matrix that is required to identify the fourth row of $T$ for the constraints defined by $\vec{\phi}_{\mu_0}$ is 
                        \begin{equation}\label{Eq:alpha_beta_system}
			  \int_{\Omega_{ref}} \hat{\phi}_{\mu_0} \hat{\phi}_*^T d \Omega \, = \, \int_{\Omega_{ref}} \phi_{A} \phi_{B} \phi_{C} \left(\hat{\phi}_{\mu_0} \hat{\phi}_*^T \right)d \Omega.
                        \end{equation}
                    This is rewritten using the transformation of bases $H$:
                        \[ \int \phi_{A}\phi_{B}\phi_{C} \left(H \hat{\phi}_* \hat{\phi}_*^T \right)d \Omega.\]
                    As the transformation $H$ has only constant elements, it may be factored out of the resulting integral
                        \begin{equation}\label{Eq:System_Matrix}
			  H \int_{\Omega_{ref}} \phi_{A}\phi_{B}\phi_{C} \left({\hat{\phi}}_* {\hat{\phi}}_*^T \right)d \Omega.
			\end{equation}
		    To see that this matrix is non-singular first notice that as $H$ is a change of basis, it is non-singular.
                    Therefore, as the product of invertible matrices is necessarily invertible, it suffices to show that the matrix defined by the integral, which we define as $\hat{M}$, is symmetric positive definite, i.e. invertible.
                    
                    The symmetry is a byproduct of the symmetry of multiplication:
                    \[
                   \hat{M}_{i,j} =  \int_{\Omega_{ref}} \phi_{A}\phi_{B}\phi_{C} \left(\hat{\phi}_{*,\,i} \hat{\phi}_{*,\,j}\right)\,d\Omega \, = \, \int_{\Omega_{ref}} \phi_{A}\phi_{B}\phi_{C} \left( \hat{\phi}_{*,\,j} \hat{\phi}_{*,\,i}\right)\,d\Omega = \hat{M}_{j,i}.\]
                    To see that $\hat{M}$ is positive definite, consider an arbitrary vector $\vec{x}$, then 
                        \[ \vec{x}^T  \hat{M} \vec{x} \, = \, \sum_{i=1}^{\text{dim}(\mathcal{T}(p-2))} \sum_{j=1}^{\text{dim}(\mathcal{T}(p-2))} \left[ x_i \int \phi_{A}\phi_{B}\phi_{C} \left( \hat{\phi}_{*,\,i} \hat{\phi}_{*,\,j}\right)\,d\Omega  \, x_j \right].\]
                    We rewrite this as 
                        \[ \int_{\Omega_{ref}}   \phi_{A}\phi_{B}\phi_{C} \left(\sum_{i=1}^{\text{dim}(\mathcal{T}(p-2))} x_i \hat{\phi}_{*,\,i}\right) \left(\sum_{j=1}^{\text{dim}(\mathcal{T}(p-2))}x_j\hat{\phi}_{*,\,j}\right)\,d\Omega  \,  .\]
                    Therefore
                        \[ \vec{x}^T \hat{M} \vec{x} \, = \,  \int_{\Omega_{ref}} \phi_{A}\phi_{B}\phi_{C} \left( \sum_{i=1}^{\text{dim}(\mathcal{T}(p-2))} x_i\hat{\phi}_{*,\,i} \right)^2 d \Omega  .\]
                    But, $\phi_{A}$, $\phi_{B}$, $\phi_{C}$, and the sum squared are all greater than or equal to zero on $\Omega_{ref}$, so 
                        \[ \vec{x}^T \hat{M} \vec{x} \geq 0. \]
                    Furthermore, because the product $\phi_{A}\phi_{B}\phi_{C}$ is nonzero on the entirety of the interior, and the sum squared term is non-negative, $\vec{x}^T \hat{M} \vec{x}$ can only be zero if $\vec{x} \, =\,0$.
                    With that, it has been shown that $\hat{M}$ is symmetric positive definite and invertible.
                    And, as the product of two invertible matrices is invertible, $H \hat{M}$ is necessarily non-singular.                
                \end{proof}
            \end{theorem}
            
            Given the invertibility of $H \hat{M}$,  we may now prove that no $(p-1)$-exact diagonal $C^0$ TFEM exist.  
            \begin{theorem}{There are no continuous triangular $(p-1)$-exact diagonal pseudo-mass matrix methods for $\mathcal{T}(p)$.}\label{Thm:p-1} 
	      \begin{proof} 
		To begin, assume that there is a basis $\vec{\psi}$ for $\mathcal{T}(p)$ for which there exists an associated $(p-1)$-exact diagonal pseudo-mass matrix.
		This basis $\vec{\psi}$ is a transformation of the modified Dubiner basis $\vec{\phi}$ under a change of basis $T$ described in \eqref{Eq:Transformation_Shape}.
		In particular, the first edge function of $\vec{\psi}$, $\psi_4$ is defined by $\vec{T}^*_4$ in \eqref{Eq:Contradiction_Constraints}.
		But Thm.~\ref{Theorem:Fourth_Row_of_T_System} shows that $\vec{T}^*_4$ is zero, which contradicts the invertibility of $T$.
	      \end{proof}
            \end{theorem}   
       \end{subsubsection}
 \end{subsection}
\end{section}

\begin{section}{A Lower-Triangular Method}

    In this section, we define a $(p-1)$-exact, lower-triangular pseudo-mass matrix approach for $p=3$.
    The intent of this section is not to develop a general method for arbitrary $p$, but rather to demonstrate that further work could lead to a $(p-1)$-exact triangular finite element method that does not require the inversion of a full mass matrix.
    That is, this section demonstrates hope for arbitrarily high-order explicit methods by extending the tools created in the contradiction of a $(p-1)$-exact diagonal pseudo-mass matrix to define other less restrictive approaches.
    
    In \cite{Helenbrook2009}, Helenbrook found a unique set of vertex functions for $\mathcal{T}(p)$ whose coefficients could be determined exactly by a diagonal operation when $f \in V_{p-1}^h$.  
    Since there is no diagonal approach to determine edge function coefficients, we relax the diagonality constraint of the pseudo-mass matrix in order to define a higher-order method that is appropriate for explicit time-stepping.  
    We rearrange the basis functions and consequent degrees of freedom to emphasize the desired lower-triangular structure.  Let $L$ be the single element pseudo-mass matrix for $p=3$ defined as  
           \begin{equation}\label{Eq:L_for_p_3}
     L = \left[ \begin{array}{ccc ccc ccc c}
			  v & 0 & 0   & 0 & 0 & 0   & 0 & 0 & 0    & 0 \\
			  0 & v & 0   & 0 & 0 & 0   & 0 & 0 & 0    & 0 \\ 
			  0 & 0 & v   & 0 & 0 & 0   & 0 & 0 & 0    & 0 \\
			  e_{1,o} & e_{1,-} & e_{1,+}   & e_1 & 0 & 0    & 0 & 0 & 0 & 0\\
			  e_{1,+} & e_{1,o} & e_{1,-}   & 0 & e_1 & 0    & 0 & 0 & 0 & 0\\
			  e_{1,-} & e_{1,+} & e_{1,o}   & 0 & 0 & e_1   & 0 & 0 & 0 & 0\\
			  e_{2,o} & e_{2,-} & e_{2,+}   & 0 & e_{2,1+} & e_{2,1-}   & e_2 & 0 & 0    & 0\\
			  e_{2,+} & e_{2,o} & e_{2,-}   & e_{2,1-} & 0 & e_{2,1+}   & 0 & e_2 & 0    & 0\\
			  e_{2,-} & e_{2,+} & e_{2,o}   & e_{2,1+} & e_{2,1-} & 0   & 0 & 0 & e_2    & 0\\
			  i_{1,A} & i_{1,B} & i_{1,C} & i_{1,{\cal A}_1} & i_{1,{\cal B}_1} & i_{1,{\cal C}_1}   & i_{1,{\cal A}_2} & i_{1,{\cal B}_2} & i_{1,{\cal C}_2}  & i_1 \\	                                    
			\end{array} \right].
    \end{equation}
    The degrees of freedom in \eqref{Eq:L_for_p_3} are organized by vertex functions $\psi_{A}$, $\psi_{B}$, $\psi_{C}$, edge functions $\psi_{{\cal A},1}$, $\psi_{{\cal B},1}$, $\psi_{{\cal C},1}$, $\psi_{{\cal A}, 2}$, $\psi_{{\cal B},2}$, $\psi_{{\cal C},2}$, and then the interior function $\psi_{{\cal I}, 1}$.
    The first three rows of $L$ correspond to the diagonal equations to determine the vertex coefficients where $v$ is the diagonal entry.  The next three rows determine the first edge function coefficients independently for each edge, but allow a dependence on the known vertex coefficients.  $e_i$ denotes the diagonal term for the $i^{th}$ edge function and the edge-vertex couplings are denoted $e_{i,o}$, $e_{i,-}$, and $e_{i,+}$ where the $o$, $-$, or $+$ indicates the opposite, clockwise, or counterclockwise vertex from the edge, respectively.  For the second function on each edge, coupling to the known first edge coefficients are allowed.  These couplings are denoted $e_{2,1 -}$ and $e_{2,1 +}$ in which the pair $2,1$ denotes the coupling from function 2 to function 1 and the $-$ or $+$ describe the relative edge location of the function 1 coefficient (clockwise or counterclockwise from the function 2 coefficient edge). For the one interior function at $p=3$, the subscripted $i$ entries are the couplings to the known vertex and edge coefficients as indicated by the subscript. 
    
    The pseudo-mass matrix is intended to work with a basis $\vec{\psi}$, which is expressed as a change of basis, $T$, of $\vec{\phi}$, the modified Dubiner basis.  We include an additional constraint that the change of basis of the edge functions (${\cal A}_{\cal{A}},  {\cal B}_{\cal{B}}$  and ${\cal C}_{\cal{C}}$ in \eqref{Eq:Transformation_Shape}) must be upper triangular.  This implies that on any edge, only the addition of higher edge functions to lower order functions is permitted.  This restriction eases the analysis of the systems defined in this section.   In order to work along with $L$, the organization of $T$ used in this section differs from that which was used previously in the paper.  Following the ordering of $L$, the change of basis for $p=3$ is defined as
    \begin{equation}\label{Eq:Permutated_Change_of_Basis_3}
            T = \left[ \begin{array}{ccc ccc ccc c}
             1&0&0  & 0 & A_{{\cal B},1} & A_{{\cal C}, 1}& 0 & A_{{\cal B}, 2} & A_{{\cal C}, 2} & A_{{\cal I}, 1} \\
             0&1&0  & B_{{\cal A},1} & 0 & B_{{\cal C}, 1}& B_{{\cal A}, 2} & 0 & B_{{\cal C}, 2} & B_{{\cal I}, 1} \\
             0&0&1  & C_{{\cal A},1} & C_{{\cal B}, 1}& 0 & C_{{\cal A}, 2} & C_{{\cal B}, 2} & 0 & C_{{\cal I}, 1} \\
             
             0&0&0  &1&0&0  &{{\cal A}}_{{\cal A},1,2}&0&0   &{{\cal A}}_{{\cal I},1,1} \\
             0&0&0  &0&1&0  &0&{{\cal B}}_{{\cal B},1,2}&0   &{{\cal B}}_{{\cal I},1,1} \\
             0&0&0  &0&0&1  &0&0&{{\cal C}}_{{\cal C},1,2}   &{{\cal C}}_{{\cal I},1,1} \\
             
             0&0&0  &0&0&0  &1&0&0   &{{\cal A}}_{{\cal I},2,1} \\
             0&0&0  &0&0&0  &0&1&0   &{{\cal B}}_{{\cal I},2,1} \\
             0&0&0  &0&0&0  &0&0&1   &{{\cal C}}_{{\cal I},2,1} \\
             
             0&0&0  &0&0&0  &0&0&0 &1
            \end{array}\right].
        \end{equation}


    \begin{subsection}{Defining \texorpdfstring{$T$}{T} and \texorpdfstring{$L$}{L} for \texorpdfstring{$p=3$}{p=3}}\label{Sec:L_and_T} 
        The matrices $L$ and $T$ are unknown.
        Being that the desired method is $(p-1)$-exact, the pseudo-mass matrix $L$ will define a coefficient vector $\vec{u}_{\psi}$ by \[\vec{u}_{\psi}^T = L^{-1} \int_{\Omega_{ref}} f \vec{\psi} d\Omega_{ref}, \] that satisfies $\vec{u}_{\psi}^T \vec{\psi} = f$ for any $f \in \mathcal{T}(p-1)$.        
        
        Following the same arguments that led to \eqref{Eq:Psi_Soln_by_Phi}, $\vec{u}_{\psi} = T^{-T}\vec{u}_{\phi}$, where $\vec{u}_\phi$ is again the coefficients of the representation of $f$ in the modified Dubiner basis.  
        Thus,
        \begin{equation}\label{Eq:Psi_System_using_Phi}
        L\left(T^{-T} \vec{u}_{\phi} \right) = T \int_{\Omega_{ref}} f \vec{\phi} d\Omega_{ref}.
        \end{equation} 
        Letting $\mu$ denote the set of indices of the Dubiner functions of degree $\le 2$, i.e. $(p-1)$-exact, then similar to \eqref{Eq:Contradiction_Constraints_2}, we have
        \begin{equation}\label{Eq:LT_System}
        L\left(T^{-1}_{i,:}\right)^T = T \int_{\Omega_{ref}} \phi_i \vec{\phi} d\Omega_{ref} \text{ for all } i \in \mu,
        \end{equation}
        in which 
        \begin{equation}\label{Eq:mu_p3}
         \mu = [1,2,3,4,5,6].
        \end{equation}
        and $T^{-1}_{\mu,:}$ can be derived from \eqref{Eq:Inverse_Transformation_Shape} as
 
\begin{multline} \label{eq:Tinv3}
(T^{-1}_{\mu,:})^T = \\
\left[\begin{array}{cccccc}
 1 & 0 & 0 & 0 & 0 & 0 \\
0 & 1 & 0 & 0 & 0 & 0 \\
0 & 0 & 1 & 0 & 0 & 0 \\
 0 & -B_{{\cal A},1} & -C_{{\cal A},1}& 1 & 0 & 0 \\
 -A_{{\cal B},1} & 0 &-C_{{\cal B},1} &0 &1 &0 \\
 -A_{{\cal C},1} & -B_{{\cal C},1} & 0 &0 &0 &1 \\
0 & {\cal A}_{{\cal A},1,2} B_{{\cal A},1} -B_{{\cal A},2} &{\cal A}_{{\cal A},1,2} C_{{\cal A},1}-C_{{\cal A},2}  & -{\cal A}_{{\cal A},1,2} & 0 &0 \\
A_{{\cal B},1}\,{\cal B}_{{\cal B},1,2}-A_{{\cal C},2} &0 &{\cal B}_{{\cal B},1,2}\,C_{{\cal B},1}-C_{{\cal B},2} &0 &-{\cal B}_{{\cal B},1,2} &0 \\
A_{{\cal C},1}\,{\cal C}_{{\cal C},1,2}-A_{{\cal C},2} & B_{{\cal C},1}\,{\cal C}_{{\cal C},1,2}-B_{{\cal C},2} & 0 &0 &0 &-{\cal C}_{{\cal C},1,2} \\
\left[ T^{-1} \right]_{A,{\cal I}} & \left[ T^{-1} \right]_{B,{\cal I}} & \left[ T^{-1} \right]_{C,{\cal I}} & -{\cal A}_{\cal A}^{-1}{\cal A}_{{\cal I},1,1} & -{\cal B}_{\cal B}^{-1}{\cal B}_{{\cal I},1,1} & -{\cal C}_{\cal C}^{-1}{\cal C}_{{\cal I},1,1} 
 \end{array}\right]
 \end{multline}

        To determine a vertex function we start with the first row of $L$.  
        The only non-zero entry of the first row of $L$ is $v = L_{1,1}$, so just information in the first row of $T^{-T}$ is necessary.  Multiplying by the first row of $L$ for each $i \in \mu$ gives 
        \begin{equation}\label{Eq:Psi_A}
         v\delta(i,1) = T_{1,:} M_{:, i} \text{ for all } i\in \mu,
        \end{equation}
         in which $\delta(\cdot, \cdot)$ denotes the Kronecker delta function.  This corresponds to 6 linear equations in 6 unknowns -- $v$ and the 5 coefficients in the first row of $T$.  These equations can be inverted which determines the transformation to define the vertex function $\psi_{A}$.  A similar procedure can be used to find the other two vertex functions.  This defines the first three rows of $L$ and $T$.

        The first function on edge ${\cal A}$ is determined using row 4 of $L$.  For rows 4-6 of $L$, only the first 6 column entries are non-zero.  These are multiplied by the first 6 rows of $T^{-T}$ for each $i \in \mu$ (shown in \eqref{eq:Tinv3}).  The entries in the first 6 rows of this matrix are all known as they have been determined by the equations for the vertex functions.  Thus, for row 4 of $L$, the left hand side only involves linear combinations of the unknowns $e_{1,o}$, $e_{1,-}$, $e_{1,+}$ and $e_1$.   The right hand side, $T_{4,:} M_{:, i}$, involves the two unknowns ${{\cal A}}_{{\cal A},1,2}$ and ${{\cal A}}_{{\cal I},1,1}$.  This again is a linear system of 6 equations in 6 unknowns.  A similar process can be followed using rows 5 and 6 to determine the first edge function on each edge.

        Using row 7 of $L$ a similar process is followed to determine the next edge function on edge ${\cal A}$.   Now, the first 9 rows of \eqref{eq:Tinv3} are all known and the left hand side becomes linear equations in the variables $e_{2,o}$, $e_{2,-}$, $e_{2,+}$, $e_{2,1+}$, $e_{2,1-}$, and $e_2$.  The right hand side only involves the unknown  ${{\cal A}}_{{\cal I},2,1}$.  There are 6 equations in 7 unknowns, so the solution is not unique.  To obtain a unique solutions, somewhat arbitrarily, $e_2$ was chosen to be 1.  Repeating the process for rows 8 and 9 of $L$ completely determines $T$.
        
       The last row of $L$ gives 6 equations in 10 unknowns which is the dimension of the space.  One solution is to  let the last row of $L$ be equal to the last row of the exact mass matrix produced by $\vec{\psi}$ which is now determined.  The equations corresponding to the last row of $L$ will then be satisfied for all $f\in \mathcal{T}(p)$.
       
       The lower triangular pseudo-mass matrix determined by this process is 
        \begin{equation}\label{Eq:L_for_p3}
            L = \left(\begin{array}{cccccccccc} \frac{1}{30} & 0 & 0 & 0 & 0 & 0 & 0 & 0 & 0 & 0\\ 0 & \frac{1}{30} & 0 & 0 & 0 & 0 & 0 & 0 & 0 & 0\\ 0 & 0 & \frac{1}{30} & 0 & 0 & 0 & 0 & 0 & 0 & 0\\ -\frac{1}{180} & \frac{1}{360} & \frac{1}{360} & \frac{1}{90} & 0 & 0 & 0 & 0 & 0 & 0\\ \frac{1}{360} & -\frac{1}{180} & \frac{1}{360} & 0 & \frac{1}{90} & 0 & 0 & 0 & 0 & 0\\ \frac{1}{360} & \frac{1}{360} & -\frac{1}{180} & 0 & 0 & \frac{1}{90} & 0 & 0 & 0 & 0\\ 0 & \frac{157}{280} & -\frac{157}{280} & 0 & \frac{1}{210} & -\frac{1}{210} & 1 & 0 & 0 & 0\\ -\frac{157}{280} & 0 & \frac{157}{280} & -\frac{1}{210} & 0 & \frac{1}{210} & 0 & 1 & 0 & 0\\ \frac{157}{280} & -\frac{157}{280} & 0 & \frac{1}{210} & -\frac{1}{210} & 0 & 0 & 0 & 1 & 0\\ -\frac{1}{2520} & -\frac{1}{2520} & -\frac{1}{2520} & \frac{1}{2520} & \frac{1}{2520} & \frac{1}{2520} & 0 & 0 & 0 & \frac{1}{1260} \end{array}\right)
        \end{equation}
        and the corresponding change of basis $T$ defining the basis $\vec\psi$, with which $L$ functions is
        \begin{equation}\label{Eq:T_for_p3}
            T = \left(\begin{array}{cccccccccc} 1 & 0 & 0 & 0 & -\frac{9}{4} & -\frac{9}{4} & 0 & -\frac{7}{12} & \frac{7}{12} & 0\\ 0 & 1 & 0 & -\frac{9}{4} & 0 & -\frac{9}{4} & \frac{7}{12} & 0 & -\frac{7}{12} & \frac{21}{4}\\ 0 & 0 & 1 & -\frac{9}{4} & -\frac{9}{4} & 0 & -\frac{7}{12} & \frac{7}{12} & 0 & \frac{21}{4}\\ 0 & 0 & 0 & 1 & 0 & 0 & 0 & 0 & 0 & -\frac{7}{2}\\ 0 & 0 & 0 & 0 & 1 & 0 & 0 & 0 & 0 & -\frac{7}{2}\\ 0 & 0 & 0 & 0 & 0 & 1 & 0 & 0 & 0 & -\frac{7}{2}\\ 0 & 0 & 0 & 0 & 0 & 0 & 1 & 0 & 0 & 0\\ 0 & 0 & 0 & 0 & 0 & 0 & 0 & 1 & 0 & 3\\ 0 & 0 & 0 & 0 & 0 & 0 & 0 & 0 & 1 & -3\\ 0 & 0 & 0 & 0 & 0 & 0 & 0 & 0 & 0 & 1 \end{array}\right)
        \end{equation}
        
        Figure~\ref{Fig:Edge_Peaks} shows the basis functions determined by $T$.  The vertex functions are localized near their corresponding vertex similar to a nodal basis.  The edge functions remain modal along the edge but become localized near the edge in the interior of the element.   
        
            \begin{figure}[H]
            \centering

            {
            \includegraphics[width=.8\textwidth]{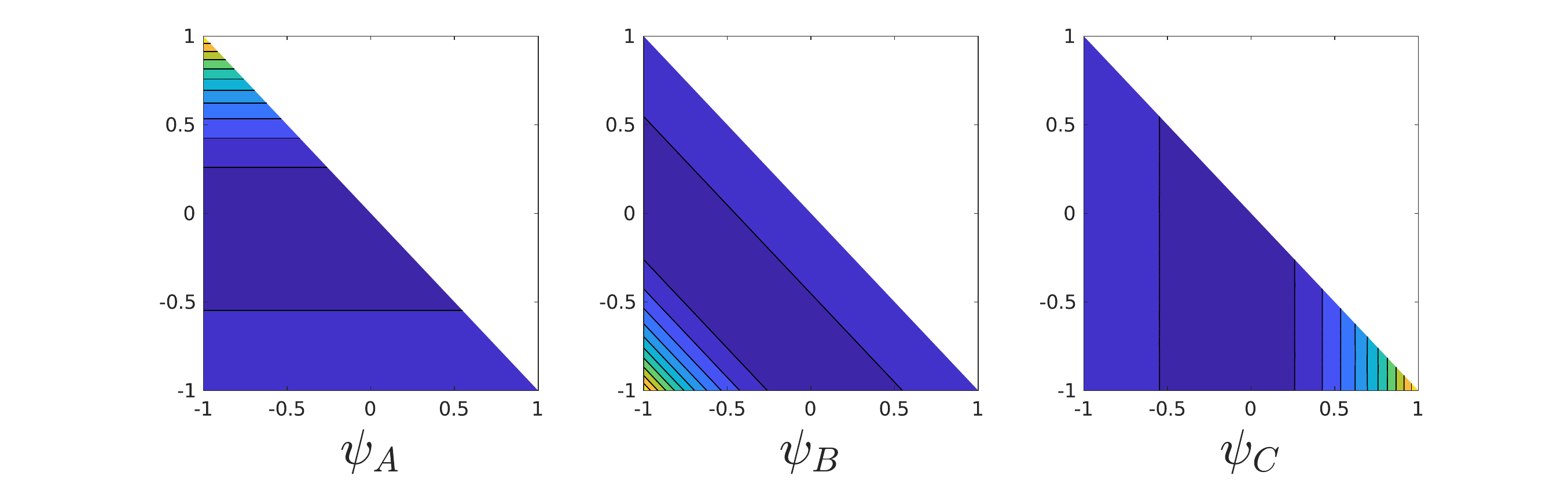}
            }
            
            {\small (a)  The vertex functions for $p=3$.}

            {
            \includegraphics[width=.8\textwidth]{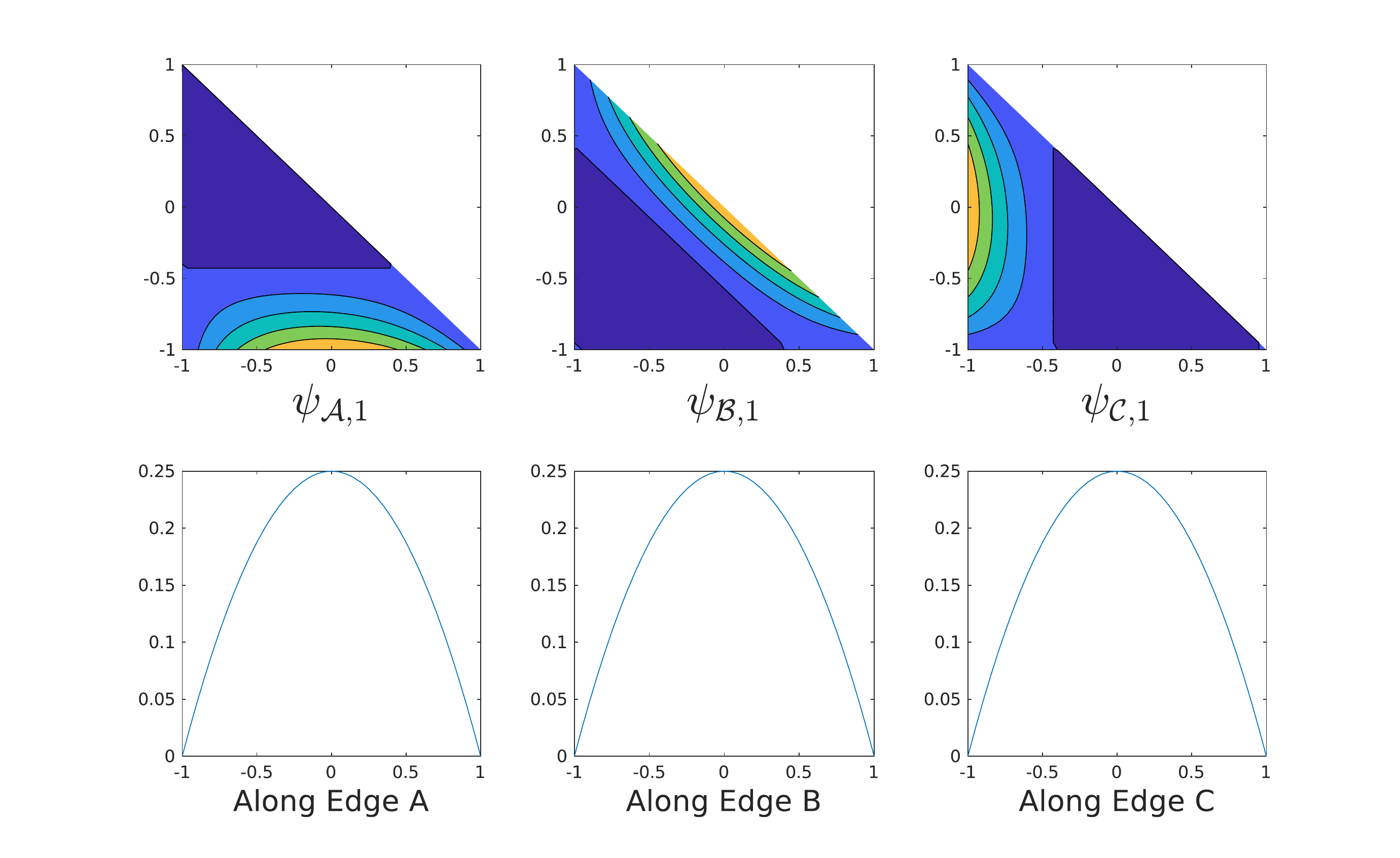}
            }
            
            {\small (b)  The quadratic edge functions for $p=3$.}
            
            {
            \includegraphics[width=.8\textwidth]{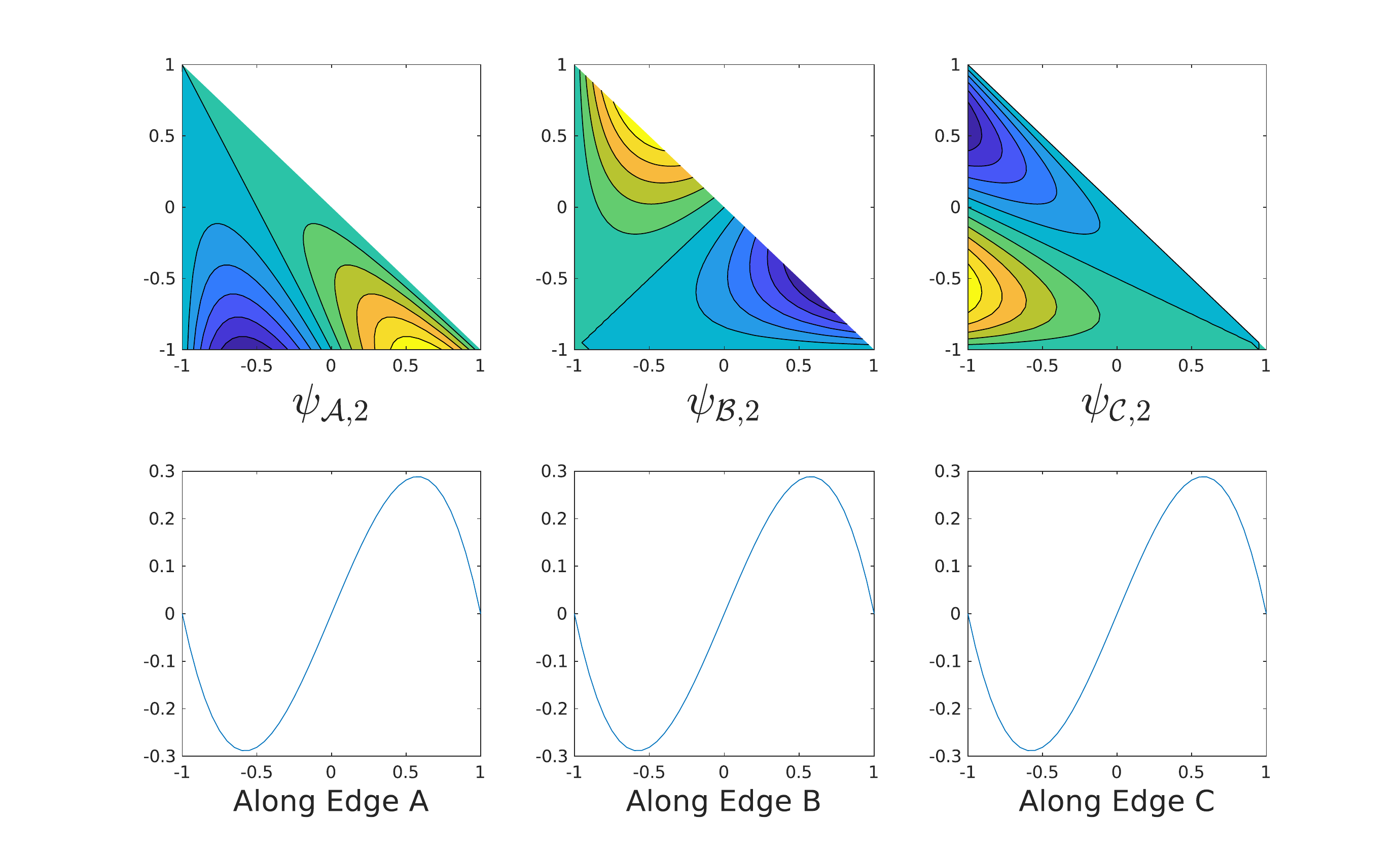}
            }
            
            {\small (c)  The cubic edge functions for $p=3$.}
            
            \caption{The basis $\vec{\psi}$ is defined as $T\vec{\phi}$, in which $\vec{\phi}$ is the modified Dubiner basis.  The contours of the vertex functions $\psi_{A}$, $\psi_{B}$, and $\psi_{C}$ in subfigure (a), demonstrate peaks at their respective vertices.  The quadratic and cubic functions $\psi_{{\cal A}, 1}$, $\psi_{{\cal B},1}$, $\psi_{{\cal C},1}$, $\psi_{{\cal A}, 2}$, $\psi_{{\cal B},2}$, and $\psi_{{\cal C},2}$ in subfigures (b) and (c) also show a rotational symmetry along the edges necessary to maintain inter-elemental continuity.  As the cubic interior function $\psi_{{\cal I}, 1}$ is $\phi_{{\cal I}, 1}$ the contour plot of this function is omitted.} \label{Fig:Edge_Peaks}
        \end{figure}
  \end{subsection}

    \begin{subsection}{A Brief Numerical Test}\label{Sec:Numerical_Analysis}
        
        This section presents a rudimentary numerical test of the $(p-1)$-exact lower-triangular method given by $L$ and $T$ in \eqref{Eq:L_for_p3} and \eqref{Eq:T_for_p3} for $p=3$.
        The method is implemented as in \eqref{Eq:Psi_System_using_Phi} via the two matrices $L$ and $T$ with the modified Dubiner basis functions.
        The method was used on a series of refined meshes to project a function outside $\mathcal{T}(3)$, and the resulting $L_2$ errors of the projection were determined.  The convergence rate of the projection is then compared to the expected rate of $3$.
                
         To implement the lower triangular approach, \eqref{Eq:Psi_System_using_Phi} is written as
         \begin{equation}\label{Eq:LT_Method}
         \vec{u}_\phi = T^T L^{-1} T \int_{\Omega} f \vec{\phi} d\Omega.
        \end{equation}
        In this case, $T$ is the global change of basis which is defined using the local $T$ matrices.  The global $L$ is assembled from element $L$ matrices in the normal FEM assembly process.   The term to the right of $L^{-1}$ is assembled element by element by integrating with respect to the Dubiner basis on each element, multiplying by the element $T$ matrix and then performing an FEM assembly.  $L$ is inverted by first inverting all of the vertex rows, which are diagonal, then doing all of the first edge functions in the mesh, then the second functions, and finally the interior functions.   Each of these operations only involves bringing the the known lower-diagonal terms of $L$ to the right hand side and then dividing by the diagonal.    This then determines the coefficients $\vec{u}_\psi$.  Multiplying by $T^T$, which again involves only local operations, determines $\vec{u}_\phi$. 
                
        The spatial convergence of the method is tested using the function
        \begin{equation}\label{Eq:Test_Function}
            f(x,y) = 1+\cos(\pi x) +\sin(\pi y)
        \end{equation}
        over a sequence of structured meshes with reducing mesh size $h$ as seen in Fig.~\ref{Fig:Mesh}.
            \begin{figure}[H]
        \centering
            \begin{tikzpicture}[thick, scale=.8]
            \draw[thick] (-1,-1)--(3,-1)--(3,3)--(-1,3)--(-1,-1);
            \foreach \i in {-1,...,2}{
            \foreach \j in {-1,...,2}{
                \draw(\i, \j) -- (\i+1, \j) -- (\i+1, \j+1) -- (\i, \j);

            }}
            \node at (-1.3,3.3) (nTL) {$(-1,1)$};
            \node at (-1.3,-1.3) (nBL) {$(-1,-1)$};
            \node at (3.3,-1.3) (nBR) {$(1,-1)$};
            \node at (3.3,3.3) (nBR) {$(1,1)$};

            \draw (0,-1.2)-- (0, -2.4); \draw (1,-1.2)--(1,-2.4);
            \draw [<-] (0,-2.2) -- (0.3, -2.2); \draw [->] (0.7,-2.2) -- (1, -2.2);
            \node at (0.5, -2.2) (nLh) {$h$};
            \draw (3.2,0)-- (4.4,0); \draw (3.2,1)-- (4.4,1);
            \draw [<-] (4.2, 0) -- (4.2, 0.3); \draw [->] (4.2, 0.7) -- (4.2, 1);
            \node at (4.2, 0.5) (nRh) {$h$};

            \draw [->](3.2, 2)--(3.8,2); \draw [->](3.3,1.9)--(3.3, 2.5);
            \node at (4,2) (nRx) {$x$}; \node at (3.3, 2.7) (nRy) {$y$};
            \end{tikzpicture}
            \caption{Structured mesh $\Omega^h$ of the spatial domain $\Omega = \left\{ (x,y) | -1 \leq x,y \leq 1\right\}$ with mesh size given as $h$.}\label{Fig:Mesh}
            \end{figure}
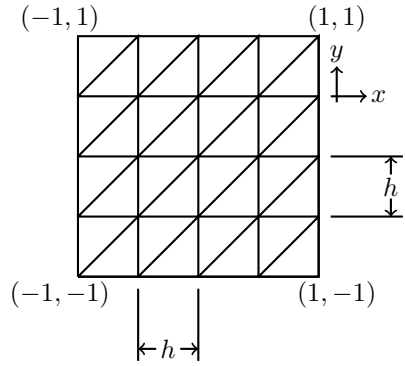
        \noindent The approximation to $f(x,y)$ is $u^h(x,y)$ as defined in \eqref{Eq:LT_Method}.
            The $L_2$ error is defined in the typical way as
            \begin{equation}\label{Eq:L2_Error}
                || f-u^h||_{L_2} = \sqrt{ \int_{\Omega} \left( f-u^h\right)^2 d\Omega }.
            \end{equation}
            This integral is numerically computed over each element of $\Omega^h$ by means of GL integration for each $h$ refinement, and the results are summarized in Table \ref{Tab:LT_Errors}.  The step size $h$ ranges from $2^{-2}$ through $2^{-7}$.  The errors are measured in the typical $L_2$-norm, \eqref{Eq:L2_Error}. The errors and step sizes are used to approximate the spatial convergence rates by $\frac{log(error_k/error_{k-1})}{log(h_k/h_{k-1})}$.  The convergence rate approaches $3$ as the meshes are refined, which is what is expected for a method that is exact for quadratic functions.
            
    \begin{table}[h]
    {\tabulinesep = 1.3mm
    \begin{tabu}{c || c | c | c | c | c | c }
     $k$ & 1 & 2 & 3 & 4 & 5 & 6 \\
     \hline
     $h$& 
         $\left( \frac{1}{2} \right)^2$ & 
        $\left( \frac{1}{2} \right)^3$ & $\left( \frac{1}{2} \right)^4$ & $\left( \frac{1}{2} \right)^5$ & 
        $\left( \frac{1}{2} \right)^6$ & $\left( \frac{1}{2} \right)^7$\\
     \hline 
     error & $4.0e-3$ & $5.3e-4$ & $6.8e-5$ & $8.6e-6$ & $1.1e-7$ & $1.3e-8$ \\
     \hline
     log-log slope & & $2.6$ & $2.8$ & $2.9$ & $3.0$ & $3.0$  
    \end{tabu}}
    \caption{Errors and convergence rates of the $p=3$ lower triangular inversion.}\label{Tab:LT_Errors}
    \end{table}

    
    \end{subsection}
\end{section}

  \begin{section}{Conclusion}\label{Sec:Conclusion}
  
    We have proven that there is no $(p-1)$-exact, diagonal pseudo-mass matrix associated with any basis $\vec{\psi}$ of $\mathcal{T}(p)$ (Thm.~\ref{Thm:p-1}).  However, by adapting the techniques used for this proof, we have developed a $(p-1)$-exact $C^0$ TFEM for $p=3$ that is appropriate for explicit time advancement as it avoids the need to perform a global inversion of the mass matrix.
    Similar to the spectral element method, this lower triangular inversion attains a $p^{th}$-order spatial convergence rate while avoiding the need to work with a large mass matrix.  In our future work, we will explore whether this method extends to arbitrarily high-order.  
    
  \end{section}

  \bibliographystyle{siamplain}
    \bibliography{references}

\end{document}